\documentclass[12pt, 14paper,reqno]{amsart}
\usepackage{amsmath,amsfonts,amssymb}

\theoremstyle{plain}
\newtheorem{theorem}{Theorem}

\newtheorem{proposition}{Proposition}
\theoremstyle{proof}
\theoremstyle{definition}

\theoremstyle{remark}

\theoremstyle{lamma}

\numberwithin{equation}{section}
\numberwithin{lemma}{section}
\numberwithin{theorem}{section}
\numberwithin{proposition}{section}

\usepackage{amsmath}
\usepackage{amsfonts}   
\usepackage{amssymb}
\usepackage{amssymb, amsmath, amsthm}
\usepackage[breaklinks]{hyperref}

\theoremstyle{thmrm}

\newenvironment{dedication}
    {\vspace{3ex}\begin{quotation}\begin{center}\begin{em}}
    {\par\end{em}\end{center}\end{quotation}}
    
\begin{document}
\title[Class groups of imaginary quadratic fields of $3$-rank $\geq 2$]{Class groups of imaginary quadratic fields of $3$-rank at least $2$}
\author{Kalyan Chakraborty and Azizul Hoque}

\address{Kalyan Chakraborty @Kalyan Chakraborty, Harish-Chandra Research Institute,
Chhatnag Road, Jhunsi, Allahabad 211 019, India.}
\email{kalyan@hri.res.in}

\address{Azizul Hoque @Azizul Hoque Harish-Chandra Research Institute,
Chhatnag Road, Jhunsi,  Allahabad 211 019, India.}
\email{ azizulhoque@hri.res.in}

\keywords{Quadratic fields, Ideal class groups, Rank}
\subjclass[2010] {Primary: 11R29, Secondary: 11R11}

\maketitle
\begin{dedication}
{\tiny To Professor Imre K\'{a}tai on the occasion of his $80^{th}$ birth anniversary, with friendship and respect.}
\end{dedication}

\begin{abstract}
We produce an infinite family of imaginary quadratic fields whose ideal class groups have $3$-rank at least $2$.
\end{abstract}

\section{Introduction}
The class group of a number field is one of the fundamental and mysterious objects in algebraic number theory. Starting from Gauss, this topic has been received serious attention of many mathematicians. It is well known that there are infinitely many imaginary quadratic fields with class number divisible by a given integer $n\geq 2$ (cf. \cite{NA22, CHKP17}). A closely related problem is concerning the $p$-rank of class groups of imaginary quadratic fields (in fact, any number fields). A result of Y. Yamamoto \cite[Proposition 2]{YA70} gives the existence of infinitely many imaginary quadratic fields whose class groups have $p$-rank at least $2$ for any integer $p\geq 2$. In \cite{DI78}, F. Diaz y Diaz developed an algorithm for generating imaginary quadratic fields whose class groups have $3$-rank at least $2$. In \cite{EKMP07}, the authors obtained a parameterized family of quadratic fields whose class group has $3$-rank at least $2$. In 2013, Y. Kishi \cite{KI13} gave a family of imaginary quadratic fields whose $3$-rank of the class group is at least 2. In \cite{LQ88}, P. Llorente and J. Quer found $1824$, $20$ and $3$ imaginary quadratic fields whose $3$-rank of the class groups are $4, 5 $ and $6$, respectively. The aim of this paper is to produce an infinitely family of imaginary quadratic fields whose class groups have $3$-rank at least $2$. 

For any three positive integers $k, \ell$ and $n$, we consider the quadratic fields:
$$
K_{-}=\mathbb{Q}(\sqrt{\ell^2-2\ell k^{3n}}) \text{ and } K_{+}=\mathbb{Q}(\sqrt{3(2\ell k^{3n}-\ell^2)}).
$$   
In this paper, we prove the following:

\begin{theorem}\label{thm1}
Let $k\equiv 4\pmod {135}$, $\ell \equiv 2\pmod {135}$ and $n$ be three odd positive integers such that $\ell<2k^{3n}$ and $\gcd(k, \ell)=1$. If $n \not\equiv 0 \pmod 3$, then the $3$-rank of the class groups of $K_{-}$ is at least $2$.
\end{theorem}
It is easy to see that Theorem \ref{thm1} yeilds infinitely many imaginary quadratic fields whose class groups has $3$-rank at least 2. 
The idea of the proof is to construct real quadratic fields of the form $K_{+}$ whose class number is divisible by $3$, and then apply the relation \cite[Theorem 1]{KI13} between the ranks of real and imaginary quadratic fields.

\section{Proof of Theorem \ref{thm1}}
We begin the proof with the following crucial proposition.
\begin{proposition}\label{proposition2.1}
Let $k, \ell$ and $n$ be as in Theorem \ref{thm1}. Then the class number of $K_{+}$ is divisible by $3$. 
\end{proposition}

The conditions $\ell<2k^{3n}$ and $n\not\equiv 0 \pmod 3$ are not necessary in Proposition \ref{proposition2.1}. Therefore, we can suppress these two conditions to get real as well as imagainary quadratic fields of the form $K_{+}$ with class number divisible by $3$. We give the proof of Proposition \ref{proposition2.1} in the most general case, that is without counting these two conditions. 

The following characterization of Y. Kishi and K. Miayke \cite[Main Theorem]{KM00} is one of the main ingredients in the proof of Proposition \ref{proposition2.1}. 
\begin{theorem}\label{thm2.2}
For any two integers $u$ and $v$, let
 \begin{equation}\label{eq2.1}
f_{u, v}(x)=x^3-uvx-u^2.
\end{equation}
If 
\begin{enumerate}\label{ctn3.1}
\item[\rm(K-1)] $u$ and $v$ are relatively prime;
\item[\rm(K-2)] $f_{u, v}(x)$ is irreducible over $\mathbb{Q}$;
\item[\rm(K-3)] discriminant $D_{f_{u, v}}$ of $f_{u, v}(x)$ is not a perfect square in $\mathbb{Z}$;
\item[\rm(K-4)] one of the following conditions holds:
\begin{enumerate}
\item[\rm (K-4.1)] $3\nmid v,$
\item[\rm (K-4.2)] $3\mid v,\hspace*{3mm} uv\not\equiv 3\pmod 9, \hspace*{3mm} u\equiv v\pm 1 \pmod 9 ,$
\item[\rm (K-4.3)] $3\mid v, \hspace*{3mm} uv\equiv 3 \pmod 9, \hspace*{3mm} u\equiv v\pm 1 \pmod { 27} ,$
\end{enumerate}
\end{enumerate}
then the normal closure of $\mathbb{\alpha}$, where $\alpha$ is a root of $f_{u,v}(x)$, is a cyclic, cubic, unramified extension of $\mathbb{K}=\mathbb{Q}(\sqrt{D_{f_{u,v}}})$; in particular, $\mathbb{K}$ has class number divisible by $3$.
Conversely, every quadratic number field $\mathbb{K}$ with class number divisible by $3$ and every unramified, cyclic and cubic extension of $\mathbb{K}$ is given by a suitable choices of integers $u$ and $v$.
\end{theorem}

\subsection*{Proof of Proposition \ref{proposition2.1}}
\noindent We choose $u=2\ell$ and $v=3k^n$. Then $\gcd(u,v)=1$ since $\gcd(k, \ell)=1$ and $\ell\equiv 2 \pmod 3$. Also by \eqref{eq2.1}, we obtain:
$$
f_{u, v}(x)=x^3-6\ell k^nx-4\ell^2.
$$
The discriminant of $f_{u,v}$ is 
$$
D_{f_{u,v}}=144\ell^2D,
$$
where $D=3\ell(2 k^{3n}-\ell).$ As both $k$ and $\ell$ are odd, we see that 
$D\equiv 3\pmod 4$, and thus $D$ is not a square in $\mathbb{Z}$.

Since $k\equiv 4\pmod 5$ and $\ell \equiv 2\pmod 5$, so that
$$f_{u, v}(x)\equiv x^3+2x-1 \pmod 5.$$
Thus $f_{u,v}(x)$ is irreducible modulo $3$ and hence it is irreducible as a polynomial with integer coefficients as well. 

We again see that $3\mid v$. As $m\equiv 4\pmod 9$, we have $m^n\equiv 1,4,7\pmod 9$ and thus $uv=6\ell m^n=3\pmod 9$. Furthermore, $$v+1=3m^n+1\equiv 4\equiv u\pmod {27}.$$
Thus we see that $f_{u,v}(x)$ satisfies the conditions (K-1)--(K-3) and (K-4.3). Therefore by Theorem \ref{thm2.2} we complete the proof of Proposition \ref{proposition2.1}. 
\qed

We now extract the following proposition from \cite[Theorem 1]{KI13} which is needed in proving Theorem \ref{thm1}.

\begin{proposition}\label{proposition2.3}
Let $d$ be a square-free positive integer such that $d\not\equiv 0\pmod 3$. Suppose $r$ and $s$ are the $3$-ranks of the class groups of $\mathbb{Q}(\sqrt{-d})$ and $\mathbb{Q}(\sqrt{3d})$, respectively. Then $r=s+1$ if and only if there does not exist a triplet $(x,y,z)\in \mathbb{Z}\times\mathbb{Z}\times\mathbb{Z}$  satisfying the following conditions:
\begin{itemize}
\item[\rm(K-5)] $x^2-4y^3=3z^2d$,
\item[\rm(K-6)] $\gcd (x,y)=1$,
\item[\rm(K-7)] $xyz\ne 0$,
\item[\rm(K-8)] $y\equiv 1 \pmod 3$ and $x^2\equiv 1, 7\pmod 9$.
\end{itemize}
\end{proposition}

\subsection*{Proof of Theorem \ref{thm1}}
Let $r$ and $s$ be the $3$-ranks of $K_{-}$ and $K_{+}$, respectively. To prove Theorem \ref{thm1}, it is sufficient to show $r=s+1$ since $s\geq 1$ by Proposition \ref{proposition2.1}.

We can express,
\begin{equation}\label{eq2.2}
\ell^2-2\ell k^{3n}=-a^2d,
\end{equation} 
where $d$ is a square-free positive integer. 

As $k\equiv 4\pmod {27}$, so that $k^{3n}\equiv 10^n\pmod {27}$. Further $10^n\equiv 10, 19\pmod{27}$ since $n\not\equiv 0\pmod 3$. Therefore by reading \eqref{eq2.2} modulo $9$, we see that $3\mid a$. Furthermore reading \eqref{eq2.2} modulo $27$, we obtain
$a^2d\equiv 9, 18\pmod {27}$ and thus $d\equiv 1,2 \pmod 3$ since $a$ is odd and $3\mid a$. 

Let us assume that $(x,y,z)\in \mathbb{Z}\times \mathbb{Z}\times \mathbb{Z}$ be such that they satisfy all the conditions (K-5)--(K-8). Then reading the condition (K-5) modulo $4$, we see that
$$
x^2\equiv 3z^2\pmod 4.
$$ 
This shows that both $x$ and $z$ are even.

Suppose that $x=2u$ and $z=2v$ for some positive integer $u$ and $v$. Then the conditions (K-5)--(K-8) imply
\begin{equation}\label{eq2.3}
u^2-y^3=3v^2d,
\end{equation}
with $\gcd(u,y)=1, ~ uyv\ne 0$ and $u^2\equiv 4,7\pmod 9$.

If $y$ is odd, then $u$ is even, and thus by \eqref{eq2.3}, we see that $v$ is odd, and hence reading \eqref{eq2.3} modulo $27$ we arrive at a contradiction. Thus $y$ is even and therefore $u$ is odd. In this case reading \eqref{eq2.3} modulo $4$, we obtain
$$
1\equiv 3v^2\pmod 4
$$
as $d\equiv 1\pmod 4$. This is not possible. Thus we complete the proof.

\section*{Acknowledgements}
The authors are grateful of Professor Carl Erickson for his interests and comments on various aspects of related works in the literature. The second author is thankful to SERB, Govt. of India for their support under N-PDF scheme (No. PDF/2017/001758).

\end{document}